\font\eusm=eusm10                   


\font\eusms=eusm7                       

\font\eusmss=eusm5                      

\font\scriptsize=cmr7

\input amstex

\documentstyle{amsppt}
  \magnification=1100
  \hsize=6.2truein
  \vsize=9.0truein
  \hoffset 0.1truein
  \parindent=2em

\NoBlackBoxes


\define\Afr{{\frak A}}                       

\define\ah{{\hat a}}                         

\define\Bfr{{\frak B}}                       

\define\Bfrbar{{\overline{\Bfr}}}            

\define\Bof{B}                               

\define\clspan{\overline\lspan}              

\define\eh{\hat e}                           

\define\eqdef{{\;\overset\text{def}\to=\;}}     

\define\Fc{{\mathchoice                      
     {\text{\eusm F}}
     {\text{\eusm F}}
     {\text{\eusms F}}
     {\text{\eusmss F}}}}

\define\Hil{{\mathchoice                     
     {\text{\eusm H}}
     {\text{\eusm H}}
     {\text{\eusms H}}
     {\text{\eusmss H}}}}

\define\Hilo{{\overset{\scriptsize o}        
     \to\Hil}}

\define\Ic{{\Cal I}}                         

\define\Icbar{{\overline{\Ic}}}              

\define\KHil{{\mathchoice                    
     {\text{\eusm K}}
     {\text{\eusm K}}
     {\text{\eusms K}}
     {\text{\eusmss K}}}}

\define\lspan{\text{\rm span}@,@,@,}         

\define\nm#1{||#1||}                         

\define\onehat{\hat 1}                       

\define\QED{$\hfill$\qed\enddemo}            

\define\restrict{\lower .3ex                 
     \hbox{\text{$|$}}}

\define\Vc{{\Cal V}}                         

\define\Vt{{\widetilde V}}                   

\define\Wc{{\Cal W}}                         

\define\xh{\hat x}                           

\define\Yc{{\Cal Y}}                         

\define\yh{\hat y}                           

\define\zt{{\tilde z}}                       


  \newcount\bibno \bibno=0
  \def\newbib#1{\advance\bibno by 1 \edef#1{\number\bibno}}
  \def\cite#1{{\rm[\bf #1\rm]}}
  \def\scite#1#2{{\rm[\bf #1\rm, #2]}}

  \newcount\notasecflag \notasecflag=0

  \newcount\secno \secno=0 \newcount\subsecno
  \def\newsec#1{\procno=0 \subsecno=0 \notasecflag=0
    \advance\secno by 1 \edef#1{\number\secno}
    \edef\currentsec{\number\secno}}
  \def\newsubsec#1{\procno=0 \advance\subsecno by 1 \edef#1{\number\subsecno}
    \edef\currentsec{\number\secno.\number\subsecno}}

  \newcount\appendixno \appendixno=0
  \def\newappendix#1{\procno=0 \notasecflag=0 \advance\appendixno by 1
    \ifnum\appendixno=1 \edef\appendixalpha{\hbox{A}}
      \else \ifnum\appendixno=2 \edef\appendixalpha{\hbox{B}} \fi
      \else \ifnum\appendixno=3 \edef\appendixalpha{\hbox{C}} \fi
      \else \ifnum\appendixno=4 \edef\appendixalpha{\hbox{D}} \fi
      \else \ifnum\appendixno=5 \edef\appendixalpha{\hbox{E}} \fi
      \else \ifnum\appendixno=6 \edef\appendixalpha{\hbox{F}} \fi
    \fi
    \edef#1{\appendixalpha}
    \edef\currentsec{\appendixalpha}}

  \newcount\procno \procno=0
  \def\newproc#1{\advance\procno by 1
   \ifnum\notasecflag=0 \edef#1{\currentsec.\number\procno}
   \else \edef#1{\number\procno}
   \fi}

  \newcount\tagno \tagno=0
  \def\newtag#1{\advance\tagno by 1 \edef#1{\number\tagno}}

\define\clo{{\Cal O}_\infty}
\define\bn{{\Bbb N}}
\define\br{{\Bbb R}}
\define\bz{{\Bbb Z}}
\define\bc{{\Bbb C}}

\redefine\Cpx{{\bc}}
\redefine\Naturals{{\bn}}
\redefine\Integers{{\bz}}
\redefine\Xc{\Yc}

\newsec{\nonapproxdiv}
 \newproc{\BarnettThm}
 \newproc{\BarnettCor}
 \newproc{\InterAlg}
 \newproc{\PINAD}
 \newproc{\Sharper}
\newsec{\crossedprod}
 \newproc{\PICond}
  \newtag{\Btensor}
 \newproc{\TensorPI}
 \newproc{\Outer}
 \newproc{\Compare}
 \newproc{\InSA}
 \newproc{\Matr}
 \newproc{\ProjUnder}
 \newproc{\IntRhoU}
 \newproc{\IntRhoUa}
\newsec{\pisimplefp}
 \newproc{\AMNB}
 \newproc{\PropertyQDef}
 \newproc{\PropertyQHolds}
  \newtag{\HilBlast}
  \newtag{\HilAlast}
  \newtag{\TxnImage}
  \newtag{\PWnSum}
 \newproc{\twotwoCompress}
 \newproc{\ABMult}
  \newtag{\abmult}
 \newproc{\ABTMult}
 \newproc{\GeneratedBy}
  \newtag{\generatedby}
 \newproc{\BfrAlg}
 \newproc{\lerhon}
 \newproc{\NotInner}
 \newproc{\IcnIdeals}
  \newtag{\IcnQuotient}
 \newproc{\WcVanish}
  \newtag{\iconeis}
  \newtag{\drestrict}
  \newtag{\zrestrict}

\newbib{\Avitzour}
\newbib{\knBKR}
\newbib{\knBlH}
\newbib{\knB}
\newbib{\knBP}
\newbib{\CuntzOn}
\newbib{\knCu}
\newbib{\knGH}
\newbib{\knKi}
\newbib{\knKis}
\newbib{\knRth}
\newbib{\knRtw}
\newbib{\knR}
\newbib{\VDNbook}
\newbib{\knZ}

\topmatter
  \title Purely Infinite Simple $C^*$-algebras Arising from Free
         Product Constructions
  \endtitle

  \author Kenneth J\. Dykema$^{\text{\dag}}$
          \thanks\dag{} Support of an NSF Postdoctoral Fellowship and a Fields
            Institute Fellowship is gratefully acknowledged.
          \endthanks 
    and Mikael R\o rdam \endauthor

  \leftheadtext{Dykema and R\o rdam}

  \rightheadtext{Purely infinite simple $C^*$--algebras}

  \address Dept.~of Mathematics and Computer Science,
           Odense University,
           DK-5230 Odense M, Denmark
  \endaddress

  \affil Department of Mathematics \\ University of California \\
         Berkeley CA 94720--3840, U.S.A\. \\ 
         (present address: same as Mikael R\o rdam's) \\
         e-mail: dykema\@imada.ou.dk \\
       and \\
         Dept.~of Mathematics and Computer Science \\
         Odense University \\
         DK-5230 Odense M, Denmark
  \endaffil

  \abstract
Examples of simple, separable, unital, purely infinite
$C^*$--algebras are constructed, including:
\roster
\item"(1)" some that are not approximately divisible;
\item"(2)" those that arise as crossed products of any of a certain class of
$C^*$--algebras by any of a certain class of non--unital endomorphisms;
\item"(3)" those that arise as reduced free products of pairs of
$C^*$--algebras with respect to any from a certain class of states.
\endroster

  \endabstract


\endtopmatter

\document \TagsOnRight \baselineskip=18pt

\newpage

\noindent{\bf Introduction}
\vskip3ex

  We construct three classes of examples of purely
infinite, simple, unital $C^*$--algebras, which may
be of special interest.  Some of these constructions use Voiculescu's
theory freeness and his construction of reduced free products of operator
algebras, (see~\cite{\VDNbook}, see also~\cite{\Avitzour}).
The first class of examples consists of separable, purely
infinite, simple, unital
$C^*$--algebras which are not approximately divisible in the sense
of~\cite{\knBKR}.
These are the first such examples, and they are
constructed by applying a theorem of L\.~Barnett~\cite{\knB} concerning
free products of von Neumann algebras, and using an
enveloping result proved in this paper.
The existence of $C^*$--algebras with these properties was claimed
in~\scite{\knBKR}{Example 4.8}, but the proposed proof was later seen to be
slightly deficient.
With Kirchberg's result, which entails that all nuclear, simple, purely
infinite $C^*$--algebras are approximately divisible, (see~\S1), these examples
have become more important and deserving of a complete and correct description.

  Skipping ahead, the third class of examples consists of reduced free products
of $C^*$--algebras.
Relatively little is understood about the order structure of the $K_0$ group of
reduced free product $C^*$--algebras, and the knowledge of whether 
projections in these $C^*$--algebras are finite or infinite is sporadic.
Thus, they are worthy objects of 
interest,
particularly in light of the open question of whether
a simple $C^*$--algebra can be infinite but not purely infinite.
We investigate a certain class of reduced free products involving non--faithful
states, namely
$$
(\Afr,\varphi)=(A,\varphi_A)*(M_n(\Cpx)\otimes B,\varphi_n\otimes\varphi_B),
$$
where $A\neq\Cpx$ and $B$ are $C^*$--algebras with states $\varphi_A$ and
$\varphi_B$,
and where $\varphi_n$ is the state on $M_n(\Cpx)$ whose support is a minimal
projection.
We show  that $\Afr$ can be realized as
the $n\times n$ matrices over
the crossed product of a $C^*$--algebra by an endomorphism.
We go on to show that, under fairly mild hypotheses, $\Afr$ is purely infinite
and simple.

  Both the second and (as mentioned above) the third class of examples involve
crossed products of unital $C^*$--algebras by non-unital endomorphisms.
These are thus in
the spirit of Cuntz's presentation of the algebras
$O_n$~\cite{\CuntzOn}.  We give sufficient conditions for such a crossed
product to be purely infinite and simple.  The second class of examples is the
case of the crossed
product, $A\rtimes_\sigma\Naturals$, of $A=\bigotimes_1^\infty B$,
for $B$ simple and unital,
by the endomorphism
$\sigma(a_1\otimes a_2\otimes\cdots)=p\otimes a_1\otimes a_2\otimes\cdots$,
for $p\in B$ a proper projection.
Indeed, the
Cuntz algebra $O_n$ is obtained when $B=M_n(\Cpx)$ and $p$ is a minimal
projection.
We show that all such $A\rtimes_\sigma\Naturals$ are purely infinite and
simple.

\vskip3ex
\noindent{\bf Acknowledgement}
\vskip1.5ex

We would like to thank the Fields Institute and the organizers of the special
year in Operator Algebras there, where our collaboration began.

\vskip3ex
\noindent{\bf\S\nonapproxdiv. Non--approximately divisible $C^*$--algebras}
\vskip3ex

We show in this section that a theorem of L.~Barnett implies that there
exist purely infinite, simple, separable, unital  $C^*$--algebras that
are not approximately divisible. In particular, if $A$ is such a
$C^*$--algebra, then $A$ is not isomorphic to $A\otimes\clo$.

E.~Kirchberg has recently proved that if $A$ is a purely infinite, simple,
separable, nuclear, unital $C^*$--algebra, then there exists a sequence
of unital $^*$--homomorphisms $\mu_n:\clo\to A$ such that
$\mu_n(b)a-a\mu_n(b)\to 0$ as $n\to\infty$ for all $a\in A$ and $b\in
\clo$. He concludes from this that $A$ is isomorphic to $A\otimes\clo$
(see \cite{\knKi}). Since $\clo$ itself is purely infinite, simple, separable,
nuclear and unital, it follows that $\clo$ is isomorphic
to $\clo\otimes\clo$. Hence $A$ is isomorphic to $A\otimes\clo$ if and
only if $A$ is (isomorphic to) $B\otimes\clo$ for some $C^*$--algebra
$B$.
In \cite{\knBKR} a unital $C^*$--algebra $A$ is called approximately
divisible if there exists a sequence (or a net, if $A$ is non--separable)
of unital $^*$--homomorphisms $\mu_n \! :M_2(\bc)\oplus M_3(\bc)\to A$ such
that $\mu_n(b)a-a\mu_n(b)\to 0$ as $n\to\infty$ for all $a\in A$ and $b\in
M_2(\bc)\oplus M_3(\bc)$. Since there is a unital embedding of
$M_2(\bc)\oplus M_3(\bc)$ into $\clo$, we conclude from the remarks
above that each $C^*$--algebra, which is isomorphic to $A\otimes\clo$ for some 
unital $C^*$--algebra $A$, is approximately divisible.

\proclaim{Theorem \BarnettThm}{\rm {(L.~Barnett \cite{\knB})}}. There is a
type III--factor ${\Cal M}$ with a faithful normal state $\varphi$ and
with elements $a,b,c\in{\Cal M}$ such that
$$
\|x-\varphi(x)\cdot 1\|_\varphi \le 14
\max\{\|[x,a]\|_\varphi , \|[x,b]\|_\varphi , \|[x,c]\|_\varphi\}
$$
for every $x\in{\Cal M}$.
\endproclaim

\proclaim{Corollary \BarnettCor}
Let ${\Cal M}$ and $a,b,c\in{\Cal M}$ be as above. Suppose $A$ is a
unital $C^*$--subalgebra of ${\Cal M}$ which contains $a$, $b$ and
$c$. Then $A$ is not approximately divisible.
\endproclaim
\demo{Proof}
Suppose, to reach a contradiction, that $A$ is approximately
divisible. Then there is a unital $^*$--homomorphism $\mu:M_2(\bc)\oplus
M_3(\bc)\to A$, such that
$$
\max\{\|[a,\mu(x)]\| , \|[b,\mu(x)]\| ,
\|[c,\mu(x)]\|\} \le \frac{1}{30} \|x\|
$$
for all $x\in M_2(\bc)\oplus M_3(\bc)$. There is a projection $e\in
M_2(\bc)\oplus M_3(\bc)$ such that $1/3\le\varphi(\mu(e))\le 1/2$, where
$\varphi$ is the faithful normal state from Theorem \BarnettThm{}. Indeed, the
set $\Gamma$ of projections $(p,q)\in M_2(\bc)\oplus M_3(\bc)$, where
$\text{dim}(p)=\text{dim}(q)=1$ is connected. Moreover, there exist
$e_1,e_2,e_3\in\Gamma$ such that $1=e_1+e_2+e_3'$ and $0\le
e_3'\le e_3$, whence $1/3\le\varphi(\mu(e))$ and $\varphi(\mu(f))\le
1/2$ for some $e,f \in \Gamma$.

Put $p=\mu(e)$. From Theorem \BarnettThm{} we get $\|p-\varphi(p)1\|_\varphi\le
14/30$ since $\|\cdot\|_\varphi\le\|\cdot\|$. On the other hand, since
$p$ is a projection, we have
$$
\|p-\varphi(p)1\|^2_\varphi = \varphi(p)-\varphi(p)^2\ge 2/9 ,
$$
a contradiction.
\QED


Here is the enveloping result mentioned in the introduction.
\proclaim{Proposition \InterAlg}
Let $B$ be a unital (non--separable) $C^*$--algebra, and let $X$ be a
countable subset of $B$.
\roster
\item"(i)" If $B$ is simple and purely infinite, then there exists a
separable, unital, simple and purely infinite $C^*$--algebra $A$ such that
$X\subseteq A\subseteq B$.
\item"(ii)" If $B$ is nuclear, then there exists a separable,
unital, nuclear $C^*$--algebra $A$ such that
$X\subseteq A\subseteq B$.
\item"(iii)" If $B$ is simple, purely infinite and nuclear, then there
exists a separable, unital, simple, purely infinite and nuclear
$C^*$--algebra $A$ such that $X\subseteq A\subseteq B$.
\endroster
\endproclaim


\demo{Proof}
The proofs of (i) and (ii) are easily obtained from the proof given
below of (iii).

Suppose $B$ is simple, purely infinite and nuclear. We may assume that
$1_B\in X$ so that $X\subseteq A$ will imply that $A$ is unital. Recall that a
unital $C^*$--algebra $D$ is simple and purely infinite if and only if
for each positive, non--zero $a\in D$ there exists $x\in D$ with
$xax^*=1$. Moreover, $x$ above can be chosen to have norm less than
$2\|a\|^{-1/2}$. 

We will show how to construct a sequence $X=X_0\subseteq X_1\subseteq
X_2\subseteq\cdots$ of countable subsets of $B$, a sequence
$A_1\subseteq A_2\subseteq A_3\subseteq\cdots$ of separable
$C^*$--subalgebras of $B$, and a sequence $\Phi_0\subseteq
\Phi_1\subseteq \Phi_2\subseteq \cdots$ of countable families of
completely positive, finite rank contractions from $B$ into $B$ such
that the following four conditions hold for every $n\ge 0$.
\roster
\item"(a)" For all $\varepsilon\ge 0$ and for each finite subset $F$ of
$X_n$ there exists $\varphi\in\Phi_n$ such that
$\|\varphi(x)-x\|<\varepsilon$ for all $x\in F$.
\item"(b)" $\varphi(X_n)\subseteq A_{n+1}$ for all $\varphi\in\Phi_n$.
\item"(c)" $X_{n+1}$ is a dense subset of $A_{n+1}$, and $X_{n+1}\cap
A_{n+1}^+$ is a dense subset of $A_{n+1}^+$.
\item"(d)" For each positive, non--zero $a\in X_n$ there exists an $x\in
A_{n+1}$ such that $\|x\|\le 2\|a\|^{-1/2}$ and $xax^*=1$.
\endroster

Indeed, given $X_n$, since there are only countably many finite subsets
of $X_n$, by the Choi--Effros characterization of nuclearity, we can find
a countable family $\Phi_n$ of completely positive finite rank
contractions satisfying (a). Moreover, if $n\ge 1$, then we may insist
that $\Phi_n\supseteq\Phi_{n-1}$. For each positive, non--zero element
$a$ in $B$ choose $x(a)\in B$ such that $\|x(a)\|\le 2\|a\|^{-1/2}$ and
$x(a)a\, x(a)^*=1$. Suppose $X_n$, $\Phi_n$ and $A_n$ are given (where
$A_0=\{0\}$). Let $A_{n+1}$ be the $C^*$--algebra generated by $A_n$ and
the countable set
$$
\{\varphi(x) \mid x\in X_n\ ,\; \varphi\in\Phi_n\} \cup
\{ x(a) \mid a\in X_n\ ,\; a\ge 0, \; a \ne 0\}.
$$
Then $A_{n+1}$ is a separable $C^*$--subalgebra of $B$, $A_n\subset
A_{n+1}$, and (b) and (d) hold. Now choose a countable subset $X_{n+1}$
of $B$ such that $X_n\subseteq X_{n+1}$ and such that (c) holds. 

Set
$$
A = \overline{\bigcup^\infty_{n=1}A_n}\ ,\quad Y =
\bigcup^\infty_{n=0}X_n\ ,\quad  \Phi = \bigcup^\infty_{n=0} \Phi_n
.
$$
Then $Y$ is a countable dense subset of $A$ and $\varphi(A)\subseteq
A$ for all $\varphi\in\Phi$. Hence, by (a), $A$ is a separable, nuclear
$C^*$--subalgebra of $B$ which contains $X \: (=X_0)$. We must also show
that $A$ is simple and purely infinite. Assume $a\in A$ is positive and
non--zero. By (c) we can find a (non--zero) positive $a'\in X_n$ for some
$n\in\bn$ such that $\|a-a'\|\le\frac{1}{5} \|a\|$. By (d) there exists
$y\in A$ such that $ya'y^*=1$ and $\|y\|\le 2\|a'\|^{-1/2}$. This
implies that
$$
\|yay^*-1\| \le \|y\|^2 \cdot \|a-a'\| < 1 .
$$
Hence $yay^*$ is invertible. Set $x=(yay^*)^{-1/2}y\in A$. Then
$xax^*=1$ as desired.
\QED

\proclaim{Theorem \PINAD}
There exist $C^*$--algebras which are separable, unital, simple and
purely infinite, but not approximately divisible.
\endproclaim

\demo{Proof}
Combine Corollary~\BarnettCor{} and Proposition \InterAlg{} (i) with
$B={\Cal M}$
and $X=\{a,b,c\}$. Recall from \cite{\knCu} that every countably
decomposable type III--factor is simple and purely infinite.
\QED

Proposition \InterAlg{} also allows us to sharpen Kirchberg's result, which was
discussed at the beginning of this section, on the
approximate divisibility of nuclear, simple, purely infinite, unital 
$C^*$--algebras. Recall for this that a (non--unital) $C^*$--algebra $A$
is approximately divisible if $A$ has an approximate unit consisting of
projections and if $pAp$ is approximately divisible (in the sense of
\cite{\knBKR}) for all projections $p$ in $A$.

\proclaim{Theorem \Sharper}{\rm (cf.~Kirchberg \cite{\knKi})}.
Every nuclear, simple, purely
infinite $C^*$--algebra is approximately divisible.
\endproclaim

\demo{Proof}
Suppose $A$ is a nuclear, simple, purely infinite $C^*$--algebra. Then
$A$ has real rank zero by \cite{\knZ} and so, by \cite{\knBP}, $A$ has
an approximate unit consisting of projections. For each (non--zero)
projection $p$ in $A$, $pAp$ is nuclear, simple and purely infinite. 

Let $F$ be a finite subset of $pAp$.
By Proposition \InterAlg{} (iii) there exists a unital, separable, nuclear,
simple, purely infinite $C^*$--subalgebra $A_0$ of $pAp$ such that
$F\subset A_0$. From Kirchberg's theorem \cite{\knKi}, $A_0$ is
isomorphic to $A_0\otimes{\Cal O}_\infty$ and (so) $A_0$ is
approximately divisible. It follows that there, for each $\varepsilon >
0$, exists a unital $^*$--homomorphism $\mu:M_2(\bc)\oplus
M_3(\bc)\to A_0$ satisfying $\|\mu(b)a-a\mu(b)\|\le\varepsilon\|b\|$ for
all $a\in F$ and $b\in M_2(\bc)\oplus M_3(\bc)$. Hence $pAp$ is
approximately divisible.
\QED

\vskip3ex
\noindent{\bf\S\crossedprod. Crossed products}
\vskip3ex

Associate to each pair consisting of a unital $C^*$--algebra $A$ and an
injective endomorphism $\sigma$ on $A$ the crossed product
$A\rtimes_\sigma\bn$, which is the universal $C^*$--algebra generated by a
copy of $A$ and an isometry $s$ 
such that $sas^* = \sigma (a)$ for all $a \in A$. The isometry $s$ 
is non--unitary if $\sigma$ is not unital.

Let $\bar{A}$ be the inductive limit of the sequence
$$
A\overset\sigma\to\rightarrow A\overset\sigma\to\rightarrow
A\overset\sigma\to\rightarrow\cdots ,
$$
and let $\mu_n:A\to\bar{A}$ be the corresponding $^*$--homomorphisms,
which satisfy $\mu_{n+1}\circ\sigma=\mu_n$ and
$$
\bar{A} = \overline{\bigcup^\infty_{n=1}\mu_n(A)} .
$$
Observe that
$\mu_n:A\to\mu_n(1)\bar{A}\mu_n(1)$ is an isomorphism if and only if
$\sigma$ is a corner endomorphism, i.e.~if $\sigma:A\to\sigma(1)A\sigma(1)$ is
an isomorphism. If $\sigma$ is not a corner endomorphism, then $A$ and
$\bar{A}$ need not be stably isomorphic.

There is an automorphism $\alpha$ on $\bar{A}$ given by
$\alpha(\mu_n(a))=\mu_n(\sigma(a)) \; (= \mu_{n-1}(a))$ for $a\in A$. The map
$\mu_n:A\to\bar{A}$ extends to an isomorphism
$\hat{\mu}_n: A\rtimes_\sigma\bn\to\mu_n(1)
(\bar{A}\rtimes_\alpha\bz)\mu_n(1)$ which satisfies $\alpha \circ
\hat{\mu}_n = \mu_n \circ \sigma$.

Summarizing results from \cite{\knKis} and \cite{\knR} we get the following
sufficient conditions to ensure that crossed products by $\bz$ and by $\bn$
are purely infinite and simple.

\proclaim{Theorem \PICond}
\roster
\item"(i)" Let $A$ be a $C^*$--algebra, $A\ne\bc$, and let $\alpha$ be an 
automorphism on $A$. Suppose that $\alpha^m$ is outer for all $m\in\bn$
and that $A$ has an approximate unit of projections $(p_n)^\infty_1$
with the property that for each $n\in\bn$ and for each non--zero
hereditary $C^*$--subalgebra $B$ of $A$ there is a projection in $B$
which is equivalent to $\alpha^m(p_n)$ for some $m\in\bz$. Then
$A\rtimes_\alpha\bz$ is purely infinite and simple.
\item"(ii)" Let $A$ be a unital $C^*$--algebra, $A\ne\bc$, and let $\sigma$ be
an injective endomorphism on $A$. Let $\alpha$ be the automorphism on
$\bar{A}$ associated with $\sigma$ as described above. Suppose that
$\alpha^m$ is outer for all $m\in\bn$, and suppose that for each non--zero
hereditary $C^*$--subalgebra $B$ of $A$ there is a projection in $B$
which is equivalent to $\sigma^m(1)$ for some $m\in\bn$. Then
$A\rtimes_\sigma\bn$ is purely infinite and simple.
\endroster
\endproclaim

\demo{Proof}
(i) By \scite{\knR}{Theorem 2.1} and its proof, the claim follows if the
conclusions of Lemmas 2.4 and 2.5 in \cite{\knR} hold. Now, \cite{\knR},
Lemma 2.4, holds whenever $\alpha^m$ is outer for all $m\in\bn$. Secondly,
the conclusion of \cite{\knR}, Lemma 2.5, is equivalent to the assertion
that every non--zero hereditary $C^*$--subalgebra of $A$ contains a
projection, which is infinite relative to the crossed product
$A\rtimes_\alpha\bz$. To prove this from the assumptions in (i), it
suffices to show that at least one of the projections in the approximate
unit $(p_n)_1^\infty$ is infinite in $A\rtimes_\alpha\bz$.

Since $A\ne\bc$ there is an $n\in\bn$ such that $p_nAp_n\ne\bc p_n$. Let
$B$ be a non--trivial hereditary $C^*$--subalgebra of $p_nAp_n$
and let $q$ be a projection in $B$ which is equivalent to
$\alpha^m(p_n)$ for an appropriate $m\in\bz$. Because $\alpha^m(p_n)$ is
equivalent to $p_n$ in $A\rtimes_\alpha\bz$ and $q$ is a proper
subprojection of $p_n$, we conclude that $p_n$ is infinite.

(ii) Since $A\rtimes_\sigma\bn$ is isomorphic to $\mu_1(1)
(\bar{A}\rtimes_\alpha\bz) \mu_1(1)$, it suffices to show that
$\bar{A}\rtimes_\alpha\bz$ is simple and purely infinite. Set
$p_{n-m}=\mu_n(\sigma^m(1))$. Then $(p_n)_{n \in \bz}$ is an
increasing approximate unit of projections for $\bar{A}$, and
$\alpha(p_n)=p_{n-1}$. It suffices, by (i), to show that each non--zero
hereditary $C^*$--subalgebra $B$ of $\bar{A}$ contains a projection
equivalent to $p_n$ for some $n\in\bz$. Equivalently, we must show that
for each non--zero positive $a$ in $\bar{A}$, we have $p_n=xax^*$ for some
$n\in\bz$ and some $x\in\bar{A}$. Find $m\in\bn$ and $b$ in $A^+$ such
that $\|\mu_m(b)-a\|\le \frac{1}{2} \|a\|$. It follows from
\cite{\knRtw}, 2.2 and 2.4, that $yay^*=\mu_m(c)$ for some non--zero
positive $c$ in $A$ and some $y$ in $\bar{A}$. By assumption,
$\sigma^k(1)=zcz^*$ for some $k\in\bn$ and some $z\in A$. Hence
$p_n=xax^*$ when $n=m-k$ and $x=\mu_m(z)y$. 
\QED

We shall consider the following more specific example. Let $B$ be a
simple, unital $C^*$--algebra which contains a non--trivial and proper
projection $p$. Set
$$
A = \bigotimes^\infty_{j=1} B ,
\tag{\Btensor} $$
and let $\sigma$ be the injective endomorphism on $A$ given by
$\sigma(a)=p\otimes a$.
In~(\Btensor), one must take tensor product norms making $A$ into a
$C^*$--algebra such that $\sigma$ exists and is injective.
This is possible, for example, by using always $\otimes_{\min}$ or always
$\otimes_{\max}$.

\proclaim{Theorem \TensorPI}
With $A$ and $\sigma$ as above, the crossed product $A\rtimes_\sigma\bn$
is simple and purely infinite.
\endproclaim

The theorem is proved in a number of lemmas that verify that the
conditions in Theorem {3.1} (ii) hold.

\proclaim{Lemma \Outer}
If $\alpha$ is the automorphism on $\bar{A}$ associated to $\sigma$, then
$\alpha^m$ is outer for every $m\in\bn$.
\endproclaim

\demo{Proof} 
If $\alpha^m$ were inner, then $\alpha^m(a)=a$ for some non--zero
$a\in\bar{A}$. As before, set $p_{n-m}=\mu_n(\sigma^m(1))$ and let 
$$
e_n = 1\otimes 1\otimes\cdots\otimes 1\otimes p\otimes 1\otimes\cdots \in
A,
$$
with $p$ in the $n$'th tensor factor. Then 
$\|p_na p_n-a\|$ tends to $0$ and $\|\mu_n(1-e_{2n})a\|$ tends to
$\|a\|$ as $n$ tends to infinity, and 
$p_{-n}$ is orthogonal to
$\mu_n(1-e_{2n})$ for all $n$.
Because
$$
\|p_{n+m} a p_{n+m}-a\|  =  \|\alpha^m(p_{n+m} a
p_{n+m}-a)\| 
  =  \|p_n a p_n-a\| ,
$$
it follows that $a=p_na p_n$ for all $n\in\bz$. Hence
$\mu_n(1-e_{2n})a=0$ for all $n\in\bn$, which entails that $a=0$, in
contradiction with our assumptions.
\QED

For the remaining part of the proof of Theorem \TensorPI{} we need to
consider comparison
theory for positive elements as described in \cite{\knCu},
\cite{\knBlH} and \cite{\knRtw}. We remind the reader of the basic theory.

Let $A$ be a $C^*$--algebra and set
$M_\infty(A)=\bigcup^\infty_{n=1}M_n(A)$. For $a,b\in M_\infty(A)^+$
write $a\precsim b$ if there is a sequence $(x_n)$ in $M_\infty(A)$
such that $x_nbx^*_n\to a$. This order relation extends the usual
Murray--von Neumann ordering of the projections in $M_\infty(A)$. If
$a\in A^+$ and if $p$ is a projection in $A$, then $p\precsim a$ if and
only if $p$ is equivalent to a projection in the hereditary
$C^*$--subalgebra $\overline{aAa}$, which again is the case if and only
if $p=xax^*$ for some $x\in A$.

Let $a,b\in M_\infty(A)^+$. Write $a\sim b$ if $a\precsim b$ and
$b\precsim a$, and let $a\oplus b$ be the element of $M_\infty(A)^+$
obtained by taking direct sum. Put
$$
S(A) = M_\infty(A)^+/\sim\ ,
$$
let $\langle a\rangle\in S(A)$ denote the equivalence class containing
$a\in M_\infty(A)^+$, set $\langle a\rangle +\langle b\rangle = \langle
a\oplus b\rangle$ and write $\langle a\rangle\le\langle b\rangle$ if
$a\precsim b$. Then $(S(A),+,\le)$ is an abelian preordered
semigroup. Let $DF(A)$ be the set of states on $S(A)$, i.e.~the set of
additive, order preserving functions $d:S(A)\to \br$ such that
$\sup\{d(\langle a\rangle)\mid a\in A\}=1$.

\proclaim{Lemma \Compare}{\rm (cf. \scite{\knGH}{Lemma 4.1}).}
If $A$ is a unital simple $C^*$--algebra,
and if $t,t'\in S(A)$ are such that $d(t) < d(t')$ for all $d\in DF(A)$,
then $nt\le nt'$ for some $n\in\bn$.
\endproclaim

\demo{Proof}
The set of dimension functions $DF(A)$ is weakly compact, which
entails that
$$
(c=) \; {\sup}\{d(t)/d(t') \mid d \in DF(A) \} < 1.
$$
Find $m,m' \in \bn$ such that $c < m'/m < 1$. Then $d(mt) < d(m't')$ for
all $d \in DF(A)$. By \cite{\knRtw}, 3.1, this implies that $kmt + u \le
km't'+u$ for some $k \in \bn$ and some $u \in S(A)$. Because 
$A$ is algebraically simple being a simple unital $C^*$--algebra, every
non--zero element of $S(A)$ is an order unit for $S(A)$. It follows that
$lk(m-m')t' \ge u$ for some $l \in \bn$. Repeated use of the inequality
$kmt + u\le km't'+u$ yields $lkmt + u \le lkm't'+u$.  Hence, if $n =
lkm$, then 
$$
nt \le lkmt+u \le lkm't'+u \le lkm't'+lk(m-m')t' = nt',
$$
as desired.
\QED 

\proclaim{Lemma \InSA}
Let $A$ be a unital, simple, infinite dimensional $C^*$--algebra. For
each non--zero $a$ in $A^+$ and for each $m \in \bn$ there is a non--zero
$b$ in $A^+$ such that $m \langle b \rangle \le \langle a \rangle$ in
$S(A)$. 
\endproclaim
\demo{Proof}
It suffices to show this in the case where $m=2$. Since $A$ is infinite
dimensional there exist two non--zero mutually orthogonal positive
elements $a_1$ and $a_2$ in the hereditary subalgebra
$\overline{aAa}$. Observe that $\langle a_1\rangle +\langle
a_2\rangle=\langle a_1+ a_2\rangle \le \langle a\rangle$. By
\cite{\knRth}, 3.4, there is a unitary $u$ in $A$ such that
$$
u\, \overline{a_1Aa_1}\, u^*\ \cap\ \overline{a_2Aa_2} \ne \{0\} .
$$
Let $b$ be a non--zero positive element in this intersection. Then
$b\in\overline{a_2Aa_2}$ and $u^*bu\in\overline{a_1Aa_1}$, whence 
$\langle b\rangle \le \langle a_1\rangle$ and $\langle b\rangle \le
\langle a_2\rangle$. This implies $2\langle b\rangle \le \langle
a\rangle$.
\QED

\proclaim{Lemma \Matr}
Let $A$ be a $C^*$--algebra, let $A_0$ be a $C^*$--subalgebra of $A$ and
suppose that $a,b\in A^+_0$ are such that $n\langle a\rangle\le n\langle
b\rangle$ in $S(A_0)$ for some $n\in\bn$. Suppose further that for some
$m\in\bn$ there is a set of matrix units $(e_{ij})_{1\le i,j\le n}$ in
$M_m(A\cap A_0')$ and that there is a projection $e\in A\cap A_0'$ such
that
$$
\hat{e} = \pmatrix e &&& 0 \\ & 0 && \\  &&\ddots & \\ 0 &&& 0
\endpmatrix \in M_m(A\cap A_0')
$$
is equivalent to $e_{11}$ in $M_m(A\cap A_0')$. Then $n\langle ae\rangle
\le m\langle b\rangle$ in $S(A)$.
\endproclaim

\demo{Proof}
Put
$$
\bar{a} = \pmatrix a &&& 0 \\ & a && \\ &&\ddots & \\ 0 &&&
a\endpmatrix\ ,\quad \bar{b} = \pmatrix b &&& 0 \\ & b && \\
&&\ddots & \\ 0 &&& b\endpmatrix \in M_n(A_0) .
$$
Then $\langle\bar{a}\rangle = n\langle a\rangle$ and
$\langle\bar{b}\rangle = n\langle b\rangle$, and so
$\langle\bar{a}\rangle \le \langle\bar{b}\rangle$ in $S(A_0)$. It
follows that $x_k\bar{b}x^*_k\to\bar{a}$ for some sequence $(x_k)$ in
$M_n(A_0)$. Let $x_k(i,j)\in A_0$ be the $(i,j)$'th entry of $x_k$, set
$$
\tilde{a} = \pmatrix a &&& 0 \\ & a && \\ &&\ddots & \\ 0 &&&
a\endpmatrix\ ,\; \tilde{b} = \pmatrix b &&& 0 \\ & b && \\
&&\ddots & \\ 0 &&& b\endpmatrix\ ,\; \tilde{x}_k(i,j) = 
\pmatrix x_k(i,j) &&& 0 \\ & x_k(i,j) && \\
&&\ddots & \\ 0 &&& x_k(i,j)\endpmatrix
$$
in $M_m(A_0)$, and set
$$
y_k = \sum^n_{i,j=1}\ \tilde{x}_k(i,j) e_{ij} \in M_m(A) .
$$
Then
$$
y_k\, \tilde{b}\, y_k^*  =  \sum_{i,j,\alpha}\ \tilde{x}_k(i,\alpha)\,
\tilde{b}\, \tilde{x}_k(j,\alpha)^* e_{ij} 
  \to  \tilde{a}\, \sum^n_{i=1} e_{ii} .
$$
Since $\tilde{a}$ commutes with $M_m(A\cap A'_0)$ and $\hat{e}$ is
equivalent to $e_{ii}$ in $M_m(A\cap A'_0)$ we get
$$
n\langle ae\rangle  =  n\langle\tilde{a}\hat{e}\rangle = \sum^n_{i=1}
\langle \tilde{a}e_{ii}\rangle 
  =  \langle\tilde{a}\sum^n_{i=1}e_{ii}\rangle \le
\langle\tilde{b}\rangle = m\langle b\rangle .
$$
\QED

\proclaim{Lemma \ProjUnder}
Let $A$ be a $C^*$--algebra, and let $(A_n)^\infty_1$ be an increasing
sequence of $C^*$--subalgebras of $A$ whose union is dense in $A$.
\roster
\item"(i)" For every non--zero positive element $a$ in $A$ there is an
$n\in\bn$ and a non--zero positive element $b$ in $A_n$ such that $\langle
b\rangle\le\langle a\rangle$ in $S(A)$.
\item"(ii)" If $p$ is a projection in $A_n$ and $b$ is a positive element in
$A_n$ such that $\langle p\rangle\le\langle b\rangle$ in $S(A)$, then
$\langle p\rangle\le\langle b\rangle$ in $S(A_m)$ for some $m\ge n$.
\endroster
\endproclaim

\demo{Proof}
(i) This follows easily from \cite{\knRtw}, 2.2 and 2.4.

(ii) If $\langle p\rangle\le\langle b\rangle$ in $S(A)$, then
$p=xbx^*$ for some $x$ in $A$ (by \cite{\knRtw}, 2.4). Find $m\ge n$ and
$y\in A_m$ such that $\|yby^*-p\| < 1/2$. Then $p\precsim yby^*$
(relative to $A_m$) by \cite{\knRtw}, 2.2, whence $\langle
p\rangle\le\langle b\rangle$ in $S(A_m)$.
\QED

Let in the following $A$ and $\sigma$ be as in Theorem \TensorPI{}. Set
$$
A_n = \big( \bigotimes^n_{j=1} B\big) \otimes \bc 1\otimes \bc
1\otimes\cdots\quad\subseteq A .
$$
Then $(A_n)^\infty_1$ is an increasing sequence of subalgebras of $A$
whose union is dense in $A$.

\proclaim{Lemma \IntRhoU}
There exists $m\in\bn$ such that for each $n\in\bn$ there exists 
$k\in\bn$ for which $n\langle\sigma^k(1)\rangle\le m\langle 1\rangle$ in
$S(A)$.
\endproclaim

\demo{Proof}
Observe first that $1-\sigma(1) \; (= (1-p) \otimes 1 \otimes \cdots)$ 
is non--zero. Accordingly, $\langle 1-\sigma(1)\rangle$ is an order unit for
$S(A)$, and so $\langle\sigma(1)\rangle \le m\langle 1-\sigma(1)\rangle$ for
some $m \in \bn$. Hence $(m+1)\langle\sigma(1)\rangle \le m\langle
1\rangle$, which again implies that
$(m+1)\langle\sigma^k(1)\rangle \le m\langle
\sigma^{k-1}(1)\rangle$ for all $k \in \bn$. We therefore have
$$
j\langle\sigma^{k-1}(1)\rangle \ge (j +1)\langle\sigma^k(1)\rangle
$$
for every $k\in\bn$ and every $j\ge m$. Thus
$$
m\langle 1\rangle  \ge  (m+1)\langle\sigma(1)\rangle\ge (m+2)
\langle\sigma^2(1)\rangle 
  \ge  (m+3) \langle\sigma^3(1)\rangle \ge\cdots\ ,
$$
from which the claim easily follows.
\QED

\proclaim{Lemma \IntRhoUa}
For each non--zero $a\in A^+$ there are integers $k,n\ge 1$ such that
$n\langle\sigma^k(1)\rangle\le n\langle a\rangle$.
\endproclaim

\demo{Proof}By Lemma \Compare{} it suffices to show that 
there exists an integer $k$ such that
$d(\langle\sigma^k(1)\rangle) <
d(\langle a\rangle)$ for all $d\in DF(A)$. Put
$$
c = \inf\{d(\langle a\rangle)\mid d\in DF(A)\} > 0 ,
$$
let $m\in\bn$ be as in Lemma \IntRhoU{} and find $n\in\bn$ such that
$m/n<c$. Then by Lemma \IntRhoU{} there is $k\in\bn$ such that
$n\langle\sigma^k(1)\rangle \le m \langle 1 \rangle$, and so 
$d(\langle\sigma^k(1)\rangle)<c\le d(\langle a\rangle)$ for all $d\in DF(A)$.
\QED

\demo{Proof of Theorem \TensorPI}
By Theorem \PICond{} (ii) and Lemma \Outer{} it suffices to show
that for each non--zero positive $a$ in $A$ there is a $k\in\bn$ such
that $\langle\sigma^k(1)\rangle\le\langle a\rangle$.

Let $m\in\bn$ be as in Lemma \IntRhoU{} and use Lemma \InSA{} to find
a non--zero positive $b$ in $A$ with $m\langle b\rangle\le\langle a\rangle$.
Use Lemma \ProjUnder{} (i) to find  $l\in\bn$ and a non--zero positive element
$b_1$ in $A_l$ with $\langle b_1\rangle\le\langle b\rangle$. According to
Lemma~\IntRhoUa{}
and the choice of $m$ there exist $k_1,k_2,n \in \bn$ such that
$n\langle\sigma^{k_1}(1)\rangle\le n \langle b_1\rangle$ and
$n\langle\sigma^{k_2}(1)\rangle \le m\langle 1\rangle$.
Use Lemma \ProjUnder{} (ii)
to find $j\ge\max\{l,k_1\}$ such that the inequality
$n\langle\sigma^{k_1}(1)\rangle\le n\langle b_1\rangle$ holds in
$S(A_{j})$. 
Let $\lambda$ be the $j'$th power of the one--sided Bernoulli--shift on
$A$. Then $\lambda$ is an endomorphism on $A$ whose image is equal to
$A\cap A'_j$. Observe that $n\langle\lambda(\sigma^{k_2}(1))\rangle \le 
m\langle 1\rangle$ in $S(A\cap A'_j)$. It follows that there exists a system of
matrix units $(e_{ij})_{1\le i,j\le n}$ in $M_m(A\cap A'_j)$ such
that $e_{11}$ is equivalent to
$$
\hat{e} = \pmatrix\lambda(\sigma^{k_2}(1)) &&& 0 \\ & 0 &&& \\
&&\ddots & \\ 0 &&& 0\endpmatrix\ \in M_m(A\cap A'_j)
$$
relative to $M_m(A\cap A'_j)$. 
Set $k=j+k_2 \; (\ge k_1 + k_2)$. Then
Lemma \Matr{} yields 
$$
\langle\sigma^k(1)\rangle \le
\langle\lambda(\sigma^{k_2}(1))\sigma^{k_1}(1)\rangle
\le n \langle\lambda(\sigma^{k_2}(1))\sigma^{k_1}(1)\rangle
\le m\langle b_1\rangle \le \langle a\rangle
$$
as desired.
\QED

\vskip3ex
\noindent{\bf\S\pisimplefp. Purely infinite simple free product
$C^*$--algebras}
\vskip3ex

In this section, we use Theorem~\PICond{} to show that certain
$C^*$--algebras arising as
reduced free products are purely infinite and simple.

For any $C^*$--algebra $A$ with state $\varphi$, denote the defining mapping
$A\rightarrow L^2(A,\varphi)$ by $a\mapsto\ah$.
For a Hilbert space, $\Hil$, we denote by $\KHil(\Hil)$ the $C^*$--algebra of
compact operators on $\Hil$.

\proclaim{Theorem \AMNB}
Let $A\neq\Cpx$ and $B$ be
$C^*$--algebras with states $\varphi_A$ and $\varphi_B$, respectively, whose
G.N.S\. representations are faithful.
Fix $N\in\{2,3,4,\ldots\}$, let $(e_{ij})_{1\le i,j\le N}$ be a system of
matrix units for $M_N(\Cpx)$ and
let $\varphi_N$ denote
the state on $M_N(\Cpx)$ such that $\varphi_N(e_{11})=1$.
Consider the reduced free product $C^*$--algebra,
$$
(\Afr,\varphi)=(A,\varphi_A)*(M_N(\Cpx)\otimes B,\varphi_N\otimes\varphi_B).
$$
If the pair $\bigl((A,\varphi_A),(B,\varphi_B)\bigr)$ has property~Q, defined
below, then $\Afr$ is simple and purely infinite.
\endproclaim

Several of the intermediate results in the proof of this theorem are valid also
without assuming property~Q; we will explicitly remark that we need
property~Q whenever this is the case.
Let us now define property~Q.
Write
$$ \align
\Hil_A=L^2(A,\varphi_A),&\quad\Hilo_A=\Hil_A\ominus\Cpx\onehat, \\
\Hil_B=L^2(B,\varphi_B),&\quad\Hilo_B=\Hil_B\ominus\Cpx\onehat,
\endalign $$
and let $\lambda_A$ and $\lambda_B$ denote the left actions of $A$ on $\Hil_A$,
respectively $B$ on $\Hil_B$.
Denote by $A*_rB$ the reduced $C^*$--algebra free product,
$$ (A*_rB,\varphi_{A*B})=(A,\varphi_A)*(B,\varphi_B). $$
Let $\Hil_{A*B}=L^2(A*_rB,\varphi_{A*B})$
and denote the usual left action
(see~\scite{\VDNbook}{\S1.5})
of $A*_rB$ on $\Hil_{A*B}$ by $\lambda_{A*B}$.
We have
$$ \Hil_{A*B}=\Cpx\onehat\oplus\bigoplus
\Sb n\ge1 \\ \vspace{.7ex} X_j\in\{A,B\} \\ \vspace{.7ex}
 X_j\neq X_{j+1} \endSb
\Hilo_{X_1}\otimes\cdots\otimes\Hilo_{X_n}. $$
Define the following subsets of $\Hil_{A*B}$:
$$ \align
\Fc&=\Hilo_A\oplus\bigoplus_{n\ge1}\Hilo_A\otimes
(\Hilo_B\otimes\Hilo_A)^{\otimes n} \\
\Fc_l&=\Fc\oplus(\Hilo_B\otimes\Fc) \\
\Fc_r&=\Fc\oplus(\Fc\otimes\Hilo_B).
\endalign $$
Identify $\Hil_B$ with the subspace,
$\Cpx\onehat\oplus\Hilo_B\subseteq\Hil_{A*B}$.
Let $V$ be the isometry from $\Hil_{A*B}\ominus\Hil_B$ onto
$\Fc_l\otimes\Hil_B$ that sends $\Hilo_{X_1}\otimes\cdots\otimes\Hilo_{X_n}$ to
$(\Hilo_{X_1}\otimes\cdots\otimes\Hilo_{X_n})\otimes\onehat$ if $X_n=A$ and to
$(\Hilo_{X_1}\otimes\cdots\otimes\Hilo_{X_{n-1}})\otimes\Hilo_B$ if $X_n=B$.

\proclaim{Definition \PropertyQDef}\rm
Let $p_B$ be the orthogonal projection from $\Hil_{A*B}$ onto $\Hil_B$.
We say that the pair
$\bigl((A,\varphi_A),(B,\varphi_B)\bigr)$ has {\it property Q} if
$$ V(1-p_B)\lambda_{A*B}(A*_rB)(1-p_B)V^*
\;\;\bigcap\;\;\bigl(\KHil(\Fc_l)\otimes\Bof(\Hil_B)\bigr)=\{0\}. $$
\endproclaim

\proclaim{Proposition \PropertyQHolds}
$\bigl((A,\varphi_A),(B,\varphi_B)\bigr)$ has property~Q if any of the
following conditions are satisfied:
\roster
\item"(i)" $(A,\varphi_A)=(C^*_r(G_1),\tau_{G_1})$ and
$(B,\varphi_B)=(C^*_r(G_2),\tau_{G_2})$ where $G_1$ and $G_2$ are nontrivial
discrete groups and $\tau_{G_\iota}$ is the canonical trace on
$C^*_r(G_\iota)$;
\item"(ii)" there are unitaries $u_A\in\Bof(\Hil_A)\cap\lambda_A(A)'$ and
$u_B\in\Bof(\Hil_B)\cap\lambda_B(B)'$ such that $u_A\onehat_A\perp\onehat_A$
and $u_B\onehat_B\perp\onehat_B$;
\item"(iii)" $B=\Cpx$ and $\lambda_A(A)\cap\KHil(\Hil_A)=\{0\}$;
\item"(iv)" $B$ is finite dimensional and
$\lambda_{A*B}(A*_rB)\cap\KHil(\Hil_{A*B})=\{0\}$.
\endroster
\endproclaim
\demo{Proof}
It is clear that (iii)$\implies$(iv)$\implies$~property~Q.
Since the right translation operators on $l^2(G_\iota)$ are unitaries commuting
with the left translation operators, if $(A,\varphi_A)$ and $(B,\varphi_B)$ are
as in condition~(i) then condition~(ii) is satisfied.
Hence we must only show that condition~(ii) implies property~Q.
Assume~(ii) is satisfied and let $0\neq x\in A*_rB$, $x\ge0$.
We will show that
$V(1-p_B)\lambda_{A*B}(x)(1-p_B)V^*\not\in\KHil(\Fc_l)\otimes\Bof(\Hil_B)$.
There is $\zeta\in\Hil_{A*B}$ such that
$\langle \lambda_{A*B}(x)\zeta,\zeta\rangle\neq0$
and we may assume without loss of generality that either
$$ \left.\aligned
&\zeta=\onehat, \\
\text{or }&\zeta\in\Hilo_B, \\
\text{or }&\zeta\in(\Hilo_A\otimes\Hilo_B)^{\otimes n},
\text{ some }n\ge1, \\
\text{or }&\zeta\in\Hilo_B\otimes(\Hilo_A\otimes\Hilo_B)^{\otimes n},
\text{ some }n\ge1,
\endaligned\qquad\right\} \tag{\HilBlast} $$
or
$$ \left.\aligned
&\zeta\in\Hilo_A, \\
\text{or }&\zeta\in(\Hilo_B\otimes\Hilo_A)^{\otimes n},
\text{ some }n\ge1, \\
\text{or }&\zeta\in\Hilo_A\otimes(\Hilo_B\otimes\Hilo_A)^{\otimes n},
\text{ some }n\ge1.
\endaligned\qquad\right\} \tag{\HilAlast} $$
Let
$$ \gather
\sigma_A:\Bof(\Hil_A)\rightarrow\Bof(\Hil_{A*B}), \\
\sigma_B:\Bof(\Hil_B)\rightarrow\Bof(\Hil_{A*B})
\endgather $$
be the ``right actions'' as in~\scite{\VDNbook}{\S1.6}, so that, by
Voiculescu's characterization of the commutant,
$\sigma_A(u_A),\,\sigma_B(u_B)\in\lambda_{A*B}(A*_rB)'\cap\Bof(\Hil_{A*B})$.
If one of the four cases in~(\HilBlast) holds then let
$w_m=(\sigma_B(u_B)\sigma_A(u_A))^m$ so that, {\it e.g\.},
$$ w_m(\Hilo_A\otimes\Hilo_B)^{\otimes n}
\subseteq(\Hilo_A\otimes\Hilo_B)^{\otimes n+m}. $$
On the other hand, if one of the three cases in~(\HilAlast) holds then let
$w_m=(\sigma_A(u_A)\sigma_B(u_B))^m$.
Then for every $m\ge1$, $Vw_m\zeta\in(P_{k_m}\Fc_l)\otimes\Hil_B$, where $P_k$
is the projection from $\Fc_l$ onto $(\Hilo_B\otimes\Hilo_A)^{\otimes k/2}$ if
$k$ is even and onto $\Hilo_A\otimes(\Hilo_B\otimes\Hilo_A)^{\otimes(k-1)/2}$
if $k$ is odd, and where $k_{j+1}=k_j+2$ for all $j\ge1$.
In particular, $w_1\zeta,\,w_2\zeta,\,w_3\zeta,\ldots$ is a sequence of
mutually orthogonal vectors, all having the same norm.
But because $w_m$ is a unitary that commutes with $\lambda_{A*B}(x)$,
$\langle\lambda_{A*B}(x)w_m\zeta,w_m\zeta\rangle
=\langle\lambda_{A*B}(x)\zeta,\zeta\rangle\neq0$
for all $m$, hence
$V(1-p_B)\lambda_{A*B}(x)(1-p_B)V^*\not\in\KHil(\Fc_l)\otimes\Bof(\Hil_B)$.
\QED

Now we prove Theorem~\AMNB{} and we begin by examining $L^2(\Afr,\varphi)$.
Note that $(\eh_{n1})_{1\le n\le N}$ is an orthonormal basis for
$L^2(M_N(\Cpx),\varphi_N)$ and that $\eh_{11}=\onehat$.
Thus
$$ \Hil_{M_N(\Cpx)\otimes B}\eqdef
L^2(M_N(\Cpx)\otimes B,\varphi_N\otimes\varphi_B)
=\bigoplus_{n=1}^N\eh_{n1}\otimes\Hil_B. $$
Let
$\Hilo_{M_N(\Cpx)\otimes B}
=\Hil_{M_N(\Cpx)\otimes B}\ominus\Cpx(\eh_{11}\otimes\onehat_B)$.
By Voiculescu's construction,
$$ \Hil\eqdef L^2(\Afr,\varphi)
=\Cpx\onehat\oplus\bigoplus
\Sb n\ge1 \\ \vspace{.7ex} X_j\in\{A,M_N(\Cpx)\otimes B\} \\ \vspace{.7ex}
  X_j\neq X_{j+1} \endSb
\Hilo_{X_1}\otimes\cdots\otimes\Hilo_{X_n}, $$
with $\Afr$ acting on $\Hil$ on the left in the usual way.
We will examine how $e_{11}\Afr e_{11}$ acts on $e_{11}\Hil$.
Identify $\Hil_B$ with
$\eh_{11}\otimes\Hil_B\subseteq\Hil_{M_n(\Cpx)\otimes B}$
and $\Hilo_B$ with $\eh_{11}\otimes\Hilo_B$ and thus identify
$\Hil_{A*B}$ with the subspace of $e_{11}\Hil\subseteq\Hil$ spanned by all the
tensors in $\Hilo_A$ and $\eh_{11}\otimes\Hilo_B$.
Consider the subspaces of $e_{11}\Hil$ defined by
$$ \align
\Vc_0&=\Hil_{A*B}\subseteq\Hil \\ \vspace{1.5ex}
\Vc_{[n]}&
=\bigoplus_{2\le k_1,\ldots,k_n\le N}\Fc_l\otimes(\Hil_B\otimes\eh_{k_11})
\otimes\Fc\otimes(\Hil_B\otimes\eh_{k_21})\otimes\cdots\otimes\Fc
\otimes(\Hil_B\otimes\eh_{k_n1}) \\ \vspace{1.5ex}
\Vc_{(n)}&
=\bigoplus_{2\le k_1,\ldots,k_n\le N}\Fc_l\otimes(\Hil_B\otimes\eh_{k_11})
\otimes\Fc\otimes(\Hil_B\otimes\eh_{k_21})\otimes\cdots\otimes\Fc
\otimes(\Hil_B\otimes\eh_{k_n1})\otimes\Fc_r.
\endalign $$
Then
$$ e_{11}\Hil=\Vc_0\bigoplus_{n\ge1}(\Vc_{[n]}\oplus\Vc_{(n)}). $$
Let $\Wc_n=\bigoplus_{k\ge n}(\Vc_{(k)}\oplus\Vc_{[k]})$ so that
$\bigcap_{n\ge1}\Wc_n=\{0\}$.

We also regard $A*_rB$ as a subalgebra of $\Afr$ in the canonical way.
Let $A_0\subseteq\ker\varphi_A$ be such that $\lspan A_0$ is norm dense in
$\ker\varphi_A$
and $\{\xh\mid x\in A_0\}$ is an orthonormal basis for $\Hilo_A$ and
let $B_0\subseteq\ker\varphi_B$ be such that $\lspan B_0$ is norm dense in
$\ker\varphi_B$
and $\{\xh\mid x\in B_0\}$ is an orthonormal basis for $\Hilo_B$.
Let
$$ \align
F&=\bigcup_{n\ge1}A_0\undersetbrace{n-1\text{ times }B_0A_0}
\to{B_0A_0\cdots B_0A_0} \\ \vspace{1.5ex}
F_l&=F\cup B_0F,
\endalign $$
so that $\{\xh\mid x\in F\}$ is an orthonormal basis for $\Fc$,
$\{\xh\mid x\in F_l\}$ is an orthonormal basis for $\Fc_l$ and
$\lspan(\{1\}\cup B_0\cup F_l\cup FB_0)$ is dense in $A*_rB$.
For $x\in F_l$ and $2\le n\le N$ let
$$ T(x,n)=e_{11}xe_{n1}\in\Afr. $$
Then $T(x,n)$ is maps $e_{11}\Hil$ onto
$$ \aligned
\Bigl(\xh\otimes&(\Hil_B\otimes\eh_{n1})\Bigr)
\quad\oplus\quad
\Bigl(\xh\otimes(\Hil_B\otimes\eh_{n1})\otimes\Fc_r\Bigr)\quad\oplus \\
&\oplus\bigoplus\Sb m\ge1\\\vspace{0.7ex}2\le k_1,\ldots,k_m\le N\endSb
\xh\otimes(\Hil_B\otimes\eh_{n1})
\otimes\Fc\otimes(\Hil_B\otimes\eh_{k_21})\otimes\cdots\otimes\Fc
\otimes(\Hil_B\otimes\eh_{k_m1}) \\ \vspace{1ex}
&\oplus\bigoplus\Sb m\ge1\\\vspace{0.7ex}2\le k_1,\ldots,k_m\le N\endSb
\xh\otimes(\Hil_B\otimes\eh_{n1})
\otimes\Fc\otimes(\Hil_B\otimes\eh_{k_21})\otimes\cdots\otimes\Fc
\otimes(\Hil_B\otimes\eh_{k_m1})\otimes\Fc_r.
\endaligned \tag{\TxnImage} $$
More specifically, taking an orthonormal basis for $e_{11}\Hil$ consisting of
tensors, each of the form  $\yh$ or $\yh\otimes\cdots$ for some $y\in F_l$, we
see that $T(x,n)$ maps each such element to
$\xh\otimes\eh_{1n}\otimes\yh(\otimes\cdots)$.
Thus $T(x,n)$ is a one--to--one mapping from an orthonormal basis for
$e_{11}\Hil$ onto and orthonormal set spanning the space given in~(\TxnImage),
thus is an isometry from $e_{11}\Hil$ onto this space.
Hence $(T(x,n))_{x\in F_l,\,2\le n\le N}$ is a family of isometries having
orthogonal ranges.
Moreover, the strong--operator limit
$$ P_{\Wc_n}=\sum\Sb x_1,\ldots,x_n\in F_l\\2\le k_1,\ldots,k_n\le N\endSb
T(x_1,k_1)\cdots T(x_n,k_n)T(x_n,k_n)^*\cdots T(x_1,k_1)^* \tag{\PWnSum} $$
is the projection from $e_{11}\Hil$ onto $\Wc_n$.

\proclaim{Lemma \twotwoCompress}
For every $x\in F_l$ and every $2\le n,m\le N$,
$$ e_{nn}xe_{mm}=0. $$
\endproclaim
\demo{Proof}
We have
$$ e_{mm}\Hil=(\Hil_B\otimes\eh_{m1})\quad\oplus\quad
(\Hil_B\otimes\eh_{m1})\otimes\Fc_r
\quad\oplus\quad\cdots, $$
where this formula continues as in~(\TxnImage), but without the $\xh$,
so $xe_{mm}\Hil=T(x,m)\Hil\subseteq e_{11}\Hil$.
\QED

\proclaim{Lemma \ABMult}
If $z_1,z_2\in A*_rB$ then
$$ e_{11}z_1e_{11}z_2e_{11}
\in e_{11}z_1z_2e_{11}
+\clspan\bigcup\Sb 2\le n\le N\\\vspace{0.7ex}x,y\in F_l\endSb T(x,n)BT(y,n)^*.
\tag{\abmult} $$
\endproclaim
\demo{Proof}
It will be enough to show~(\abmult) for $z_1$ and $z_2$ in a set that
densely spans $A*_rB$, hence we assume without loss of generality that
$z_j\in B\cup F_lB$.
If either $z_1\in B$ or $z_2\in B$ then
$e_{11}z_1e_{11}z_2e_{11}=e_{11}z_1z_2e_{11}$.
If $z_j=f_jb_j$ for
$f_j\in F_l$ and $b_j\in B$, ($j=1,2$), then
$$ \align
e_{11}z_1e_{11}z_2^*e_{11}
&=e_{11}z_1z_2^*e_{11}-\sum_{n=2}^Ne_{11}z_1e_{n1}e_{1n}z_2^*e_{11} \\
&=e_{11}z_1z_2^*e_{11}-\sum_{n=2}^Ne_{11}f_1e_{n1}b_1b_2^*e_{1n}f_2^*e_{11} \\
&=e_{11}z_1z_2^*e_{11}-\sum_{n=2}^NT(f_1,n)b_1b_2^*T(f_2,n)^*.
\endalign $$
\QED

\proclaim{Lemma \ABTMult}
For every $x_0\in F_l$, and $2\le n\le N$,
$$ e_{11}(A*_rB)e_{11}T(x_0,n)
\subseteq\clspan\bigcup_{x\in F_l}T(x,n)B. $$
\endproclaim
\demo{Proof}
Let $z\in A*_rB$.
Then
$$ e_{11}ze_{11}T(x_0,n)=e_{11}ze_{11}x_0e_{n1}
=e_{11}zx_0e_{n1}-\sum_{k=2}^Ne_{11}ze_{kk}x_0e_{n1}=e_{11}zx_0e_{n1}, $$
where the last equality follows from Lemma~\twotwoCompress.
Thus
$$ e_{11}ze_{11}T(x_0,n)\in e_{11}(A*_rB)e_{n1}
\subseteq\clspan\bigcup_{x\in F_l}T(x,n)B. $$
\QED

\proclaim{Lemma \GeneratedBy}
$$ \aligned
e_{11}\Afr e_{11}&=C^*(e_{11}Ae_{11}\cup e_{11}B\cup\bigcup_{2\le n\le N}
e_{11}\ker\varphi_Ae_{n1}) \\ \vspace{1.5ex}
&=C^*(e_{11}(A*_r B)e_{11}\cup\bigcup_{2\le n\le N} e_{11}F_le_{n1}) \\
=\clspan&\biggl(\bigcup\Sb p,q\ge0\\x_j,y_j\in F_l\\2\le n_j,m_j\le N\endSb
T(x_1,n_1)\cdots T(x_p,n_p)e_{11}(A*_rB)e_{11}
T(y_q,m_q)^*\cdots T(y_1,m_1)^*\biggr).
\endaligned \tag{\generatedby} $$
\endproclaim
\demo{Proof}
Because $\Afr$ is generated by $A$ and $B$ together with the system of matrix
units $(e_{ij})_{1\le i,j\le N}$, we have that
$$ e_{11}\Afr e_{11}=C^*(\bigcup_{1\le i,j\le N}e_{1i}Ae_{j1}
\cup\bigcup_{1\le i,j\le N}e_{1i}Be_{j1}). $$
Using Lemma~\twotwoCompress{} and that $B$ commutes with $e_{ij}$, the first
equality of~(\generatedby) clearly holds and then the second equality
is also clear.
Now the third equality follows from Lemmas~\ABMult{} and~\ABTMult{} and the
fact that, for $x,y\in F_l$ and $2\le n,m\le N$,
$$ T(y,m)^*T(x,n)=\cases e_{11}&\text{ if $x=y$ and }n=m \\0&\text{ otherwise.}
\endcases $$
\QED

For $p\ge1$ let
$$ \align
\Bfr_p=\clspan\biggl(&\bigcup\Sb0\le k\le p-1\\x_j,y_j\in F_l\\
2\le n_j,m_j\le N\endSb
T(x_1,n_1)\cdots T(x_k,n_k)e_{11}(A*_r B)e_{11}
T(y_k,m_k)^*\cdots T(y_1,m_1)^* \\
&\cup\bigcup\Sb x_j,y_j\in F_l\\2\le n_j,m_j\le N\endSb
T(x_1,n_1)\cdots T(x_p,n_p)BT(y_p,m_p)^*\cdots T(y_1,m_1)^*\biggr)
\endalign $$
and let $\Bfr=\overline{\bigcup_{p\ge1}\Bfr_p}$.
Hence $\Wc_n$ is an invariant subspace of $\Bfr$, for all $n\ge1$.
Based on Lemmas~\ABMult{} and~\ABTMult, we have
\proclaim{Observation \BfrAlg}\rm
Each $\Bfr_p$ and also $\Bfr$ is a $C^*$--subalgebra of $e_{11}\Afr e_{11}$.
For any fixed $x_0\in F_l$, $e_{11}\Afr e_{11}=C^*(\Bfr\cup\{T(x_0,2)\})$ and
$z\mapsto T(x_0,2)zT(x_0,2)^*$ is an endomorphism, call it $\sigma$, of $\Bfr$.
Hence $e_{11}\Afr e_{11}$ is a quotient of the universal crossed product,
$\Bfr\rtimes_\sigma\Naturals$, of $\Bfr$ by the endomorphism $\sigma$.
Therefore, once we have proved the following two propositions,
Theorem~\PICond{} will
imply that this universal crossed product is simple and purely infinite and
Theorem~\AMNB{} will be proved.
\endproclaim
\proclaim{Proposition \lerhon}
If $\bigl((A,\varphi_A),(B,\varphi_B)\bigr)$ has property~Q then
for every $z\in\Bfr$ such that $z\ge0$ and $z\neq0$,
there are $n\in\Naturals$ and $x\in\Bfr$ such that $xzx^*=\sigma^n(e_{11})$,
i.e\. $\sigma^n(e_{11})\lesssim z$.
\endproclaim
\proclaim{Proposition \NotInner}
Let $(\Bfrbar,\alpha)$ be the $C^*$--dynamical system associated to
$(\Bfr,\sigma)$ as described at the beginning of~\S\crossedprod.
For no $m\ge1$ is the automorphism $\alpha^m$ of $\Bfrbar$ inner.
\endproclaim
In order to prove these propositions, let us define, for $n\ge1$,
$$ \Ic_n=\clspan\biggl(\bigcup\Sb k\ge n\\x_j,y_j\in F_l\\2\le n_j,m_j\le N
\endSb
T(x_1,n_1)\cdots T(x_k,n_k)e_{11}(A*_rB)e_{11}
T(y_k,m_k)^*\cdots T(y_1,m_1)^*\biggr). 
$$
Clearly, $\Ic_n$ is an ideal of $\Bfr$, $\Ic_n\supseteq\Ic_{n+1}$, 
$\bigcap_{n\ge1}\Ic_n=\{0\}$ and $\Ic_n$ vanishes on $\Hil\ominus\Wc_n$.
The following lemma may be of interest, although it is not used in the sequel.
\proclaim{Lemma \IcnIdeals}
$\Bfr/\Ic_1\cong A*_rB$ and
$$ \Ic_n/\Ic_{n+1}\cong(A*_rB)\otimes\KHil\quad\text{for all }n\ge1,
\tag{\IcnQuotient} $$
where $\KHil$ appearing in~(\IcnQuotient) is the algebra of compact operators
on the infinite dimensional Hilbert space
$((\Fc_l\otimes\Cpx^{N-1})^{\otimes n})$.
\endproclaim
\demo{Proof}
Let us first consider the case of $\Bfr/\Ic_1$.
Let $\iota:A*_rB\rightarrow\Afr$ be the canonical (functorial) inclusion.
Consider the mapping $\pi:A*_rB\rightarrow\Bfr/\Ic_1$ given by
$\pi(z)=[e_{11}\iota(z)e_{11}]\in\Bfr/\Ic_1$.
By Lemma~\ABMult, $\pi$ is a $*$--homomorphism.
But since $\Ic_1$ vanishes on $\Vc_0=\Hil_{A*B}$, the map from $\Bfr/\Ic_1$
to $A*_rB$ given by $[y]\mapsto y\restrict_{\Vc_0}$ is well--defined and is
the inverse of $\pi$.
Hence $\pi$ is an isomorphism.

For $n\ge1$,
$$ \{T(x_1,k_1)\cdots T(x_n,k_n)T(y_n,l_n)^*\cdots T(y_1,l_1)^*
\mid x_j,y_j\in F_l,\,2\le k_j,l_j\le N\} $$
is a system of matrix units whose closed linear span is a copy of
$\KHil((\Fc_l\otimes\Cpx^{N-1})^{\otimes n})$ and there is an isomorphism
$\Ic_n\rightarrow\Bfr\otimes \KHil((\Fc_l\otimes\Cpx^{N-1})^{\otimes n})$
given by, for $z\in A*_rB$ and
$p\ge n$,
$$ \align
T(x_1,k_1)\cdots T(x_p,k_p)zT(y_p,l_p)^*\cdots T&(y_1,l_1)^*\mapsto \\
\mapsto
\bigl(
T(x_{n+1},k_{n+1})\cdots T(x_p,k_p)e_{11}&ze_{11}
T(y_p,l_p)^*\cdots T(y_{n+1},l_{n+1})^*\otimes \\
&\otimes T(x_1,k_1)\cdots T(x_n,k_n)T(y_n,l_n)^*\cdots T(y_1,l_1)^*\bigr).
\endalign $$
This isomorphism sends $\Ic_{n+1}$ onto
$\Ic_1\otimes\KHil((\Fc_l\otimes\Cpx^{N-1})^{\otimes n})$, so
$$\Ic_n/\Ic_{n+1}\cong(\Bfr/\Ic_1)\otimes
\KHil((\Fc_l\otimes\Cpx^{N-1})^{\otimes n})$$.
\QED

\proclaim{Lemma \WcVanish}
If $\bigl((A,\varphi_A),(B,\varphi_B)\bigr)$ has property~Q,
$z\in\Bfr$ and $z$ vanishes on $\Wc_n$ for some $n$, then $z=0$.
Consequently $\nm{P_{\Wc_n}zP_{\Wc_n}}=\nm z$ for all $z\in\Bfr$
and for all $n\ge1$.
\endproclaim
\demo{Proof}
First consider the case $n=1$.
We have $e_{11}\Hil=\Vc_0\oplus\Wc_1$ and there is a unitary
$$ U_1:\Wc_1\rightarrow\Fc_l\otimes\Cpx^{N-1}\otimes\Hil_B\otimes\Xc $$
for an infinite dimensional Hilbert space, $\Xc$, such that for $x,y\in F_l$
and $2\le n,m\le N$,
$$ T(x,n)T(y,m)^*
=0_{\Vc_0}\oplus U_1^*(f_{x,y}\otimes e_{n,m}\otimes1_{\Hil_B}\otimes1_\Xc)U_1
$$
where $f_{x,y}$ is the rank--one operator on $\Fc_l$ sending $\yh$ to $\xh$ and
where $e_{n,m}$ is the rank--one operator on $\Cpx^{N-1}$ sending the $(m-1)$st
standard basis vector to the $(n-1)$st.
Hence we have that
$$ \Ic_1\subseteq0_{\Vc_0}\oplus U_1^*\bigl(\KHil(\Fc_l)\otimes\Bof(\Cpx^{N-1})
\otimes\Bof(\Hil_B)\otimes\Bof(\Xc)\bigr)U_1. \tag{\iconeis} $$
For $d\in A*_rB$ we have $e_{11}de_{11}\restrict_{\Vc_0}=\lambda_{A*B}(d)$
under the identification of $\Vc_0$ with $\Hil_{A*B}$.
Recall (from just before Definition~\PropertyQDef) the unitary
$$ V:\Hil_{A*B}\ominus\Hil_B\rightarrow\Fc_l\otimes\Hil_B. $$
Let
$$ \Vt:(\Hil_{A*B}\ominus\Hil_B)\otimes\Cpx^{N-1}
\rightarrow\Fc_l\otimes\Cpx^{N-1}\otimes\Hil_B $$
be such that
$\Vt^*(\xi_1\otimes\xi_2\otimes\xi_3)=V^*(\xi_1\otimes\xi_3)\otimes\xi_2$,
{\it i.e\.} $\Vt$ just pushes $\Cpx^{N-1}$ through and acts like $V$ on the
rest.
Then it is easily seen for $d\in A*_rB$ that
$$ U_1\bigl(e_{11}de_{11}\restrict_{\Wc_1})U_1^*
=\Vt\bigl((1-p_B)\lambda_{A*B}(d)(1-p_B)\otimes1_{\Cpx^{N-1}}\bigr)\Vt^*
\otimes1_\Xc. \tag{\drestrict} $$
{}From the proof of Lemma~\IcnIdeals, there is a $*$--homomorphism
$\psi:\Bfr\rightarrow A*_rB$, whose kernel is $\Ic_1$,
given by $\psi(z)=z\restrict_{\Vc_0}$ and
the mapping $A*_rB\ni d\mapsto e_{11}de_{11}$ is a right inverse for $\psi$.
Hence $z=e_{11}\psi(z)e_{11}+z_1$ where $z_1\in\Ic_1$ and so, using~(\iconeis)
and~(\drestrict),
$$ U_1(z\restrict_{\Wc_1})U_1^*
\in\Bigl(\Vt\bigl((1-p_B)\lambda_{A*B}(\psi(z))(1-p_B)
\otimes1_{\Cpx^{N-1}}\bigr)\Vt^*
+\KHil(\Fc_l)\otimes\Bof(\Cpx^{N-1})\otimes\Bof(\Hil_B)\Bigr)\otimes\Bof(\Xc).
\tag{\zrestrict} $$
But assuming that property~Q holds, using~(\zrestrict) the supposition that
$z\restrict_{\Wc_1}=0$ implies that $\psi(z)=0$.
But then $z\in\Ic_1$ and $\Ic_1$ has support equal to $\Wc_1$, so $z=0$.
Hence the first part of the lemma is proved in the case $n=1$.

Now we will show, for $n\ge1$, and $z\in\Bfr$ that $z\restrict_{\Wc_{n+1}}=0$
implies $z\restrict_{\Wc_n}=0$, which when combined with the case proved above
will show that $z=0$.
For every $x_1,\ldots,x_n,y_1,\ldots,y_n\in F_l$ and
$2\le k_1,\ldots,k_n,l_1,\ldots,l_n\le N$, 
$$ T(x_n,k_n)^*\cdots T(k_1,k_1)^*zT(y_1,l_1)\cdots T(y_n,l_n)\in\Bfr $$
vanishes on $\Wc_1$, hence is equal to zero.
Therefore by~(\PWnSum){}
$P_{\Wc_n}zP_{\Wc_n}=0$ and, since $\Wc_n$ is invariant under $\Bfr$,
$zP_{\Wc_n}=0$, as required.

Now since $P_{\Wc_n}$ commutes with $\Bfr$, the mapping
$z\mapsto P_{\Wc_n}zP_{\Wc_n}$ is a $*$--homomorphism with zero kernel,
proving that $\nm{P_{\Wc_n}zP_{\Wc_n}}=\nm z$.
\QED

\demo{Proof of Proposition~\lerhon}
First suppose $0\ne z\ge0$, $z\in\Bfr_p$, some $p\ge1$.
By Lemma~\WcVanish, $z\restrict_{\Wc_p}\ne0$.
Looking at the definition of $\Wc_p$, we see that there is a canonical unitary
$$ U_p:\Wc_p\rightarrow(\Fc_l\otimes\Cpx^{N-1})^{\otimes p}
\otimes\Hil_B\otimes\Xc_p $$
for $\Xc_p$ the infinite dimensional Hilbert space,
$\Cpx\xi\oplus\Fc_r\oplus\cdots$,
such that
$$ U_p\Bfr_pU_p^*
=\KHil((\Fc_l\otimes\Cpx^{N-1})^{\otimes p})
\otimes\lambda_B(B)\otimes1_{\Xc_p} $$
with
$$ (U_pT(x_1,k_1)\cdots T(x_n,k_n)T(y_n,l_n)^*\cdots T(y_1,l_1)^*U_p^*)_{
x_j,y_j\in F_l,\,2\le k_j,l_j\le N} $$
being a system of matrix units for
$\KHil((\Fc_l\otimes\Cpx^{N-1})^{\otimes p})
\otimes1_{\Hil_B}\otimes1_{\Xc_p}$.
Since $P_{\Wc_p}$ commutes with $z$ we have
$z\ge z\restrict_{\Wc_p}=P_{\Wc_p}zP_{\Wc_p}\ge0$ and there is a vector,
$\zeta\in(\Fc_l\otimes\Cpx^{N-1})^{\otimes p}\otimes\Hil_B$ such that if $q$
is the projection onto $\zeta\otimes\Xc_p$ then $z\ge cU_p^*qU_p$ for some
$c>0$.
Let $\zeta'\in\Fc\otimes\Cpx^{N-1}$ be nonzero, so that
$\zeta\otimes\zeta'\in(\Fc_l\otimes\Cpx^{N-1})^{\otimes p+1}$
and let $q'$ be the projection onto
$\zeta\otimes\zeta'\otimes\Hil_B\otimes\Xc_{p+1}$.
Then $z\ge cU_p^*qU_p^*\ge cU_{p+1}^*q'U_{p+1}$.
Since both $q'$ and $U_{p+1}T(x_0,2)^{p+1}(T(x_0,2)^*)^{p+1}U_{p+1}^*$ are
minimal projections in
$\KHil((\Fc_l\otimes\Cpx^{N-1})^{\otimes p})
\otimes1_{\Hil_B}\otimes1_{\Xc_p}\subseteq U_{p+1}\Bfr_{p+1}U_{p+1}^*$, it is
clear that
$U_{p+1}^*q'U_{p+1}$ is equivalent in $\Bfr_{p+1}$ to the projection
$T(x_0,2)^{p+1}(T(x_0,2)^*)^{p+1}=\sigma^{p+1}(e_{11})$.
Hence $z\gtrsim\sigma^{p+1}(e_{11})$, as required, and there is
$y\in\Bfr$ such that $\sigma^{p+1}(e_{11})=yzy^*$.
Since
we may assume that $U_p^*\zeta U_p$ is in the range of
the spectral projection for $z$
corresponding to the interval $\bigl[\frac12\nm z,\nm z\bigr]$,
we may assume $\nm y\le\sqrt2\nm z^{-1/2}$.

For general $z\in\Bfr$, $0\ne z\ge0$, there is $p\ge1$ and $\zt\in\Bfr_p$,
$0\ne\zt\ge0$ such that $\nm{z-\zt}<\nm z/3$.
Let $y$ be such that $\nm y\le\sqrt2\nm\zt^{-1/2}$ and
$\sigma^{p+1}(e_{11})=y\zt y^*$ holds.
Then
$$ \nm{yzy^*-y\zt y^*}\le\frac{2\nm{z-\zt}}{\nm\zt}\le\frac{\nm z}{3\nm\zt}<1,
$$
and by the continuous function calculus and
the theory of projections there is $y'\in\Bfr$ such that
$\sigma^{p+1}(e_{11})=y'z(y')^*$, namely $\sigma^{p+1}(e_{11})\lesssim z$.
\QED

\demo{Proof of Proposition \NotInner}
Recall that $\Bfrbar$ is the inductive limit of
$$ \Bfr\overset\sigma\to\rightarrow\Bfr\overset\sigma\to\rightarrow\cdots $$
with corresponding embeddings $\mu_n:\Bfr\rightarrow\Bfrbar$ ($n\ge1$) such
that $\mu_n(z)=\mu_{n+1}(\sigma(z))$ and the
automorphism $\alpha$ is defined by $\alpha(\mu_n(z))=\mu_n(\sigma(z))$
for all $n\ge1$ and $z\in\Bfr$.
For $k\in\Integers$ consider
$$ \Icbar_k=\overline{\bigcup_{n\ge\max(1,-k)}\mu_n(\Ic_{k+n})}
\subseteq\Bfrbar. $$
Then $\Icbar_k$ is an ideal of $\Bfrbar$ and
$$ \Icbar_k\supseteq\Icbar_{k+1},\quad
\bigcup_k\Icbar_k=\Bfrbar,\quad\bigcap_k\Icbar_k=\{0\}\quad\text{and}\quad
\alpha(\Icbar_k)=\Icbar_{k+1}. $$
Let
$\pi_k:\Bfrbar\rightarrow\Bfrbar/\Icbar_k$
be the quotient map $z\mapsto z+\Icbar_k$.
Since the union of the ideals is $\Bfrbar$ and their intersection is $\{0\}$,
we have for every $z\in\Bfrbar$ that
$$ \lim_{k\rightarrow-\infty}\nm{\pi_k(z)}=0\quad\text{and}\quad
\lim_{k\rightarrow+\infty}\nm{\pi_k(z)}=\nm z. $$
Moreover,
$$ \nm{\pi_k(z)}=\nm{z+\Ic_k}=\nm{\alpha(z+\Ic_k)}
=\nm{\alpha(z)+\Icbar_{k+1}}=\nm{\pi_{k+1}(\alpha(z))}. $$
Supposing to obtain a contradiction that $\alpha^m$ is inner for some $m\ge1$,
there is
$0\neq u\in\Bfrbar$ such that $\alpha^m(u)=u$.
There are $\epsilon>0$ and $K_1,K_2\in\Integers$ such that
$\nm{\pi_k(u)}<\epsilon$
for all $k\le K_1$ and $\nm{\pi_k(u)}>2\epsilon$ for all $k\ge K_2$.
But $\nm{\pi_k(u)}=\nm{\pi_{k+m}(\alpha^m(u))}=\nm{\pi_{k+m}(u)}$, which
implies $\epsilon>2\epsilon$.
\QED


\enddocument